\documentclass[reqno, 10pt]{amsart}

\usepackage{amsfonts}

\newcommand {\ZZ}{Z \!\!\!Z}

\newcommand{\ra}{\rightarrow}
\newcommand{\no}{\noindent}
\newcommand{\oo}{\infty}
\newcommand{\nn}{\nonumber}

\newtheorem{thm}{Theorem}[section]
\newtheorem{prop}[thm]{Proposition}
\newtheorem{lemma}[thm]{Lemma}

\newtheorem{remark}{Remark}

\title{Tightness for the interfaces of one-dimensional voter models}
\author{S. Belhaouari \qquad T. Mountford \qquad G. Valle}
\address{D\'epartement de Math\'ematiques, \'Ecole Polytechnique
F\'ed\'erale, 1015 Lausanne, Switzerland. 
\newline
e-mail:  \rm \texttt{thomas.mountford@epfl.ch}
\newline
e-mail:  \rm \texttt{samir.brahim@epfl.ch}
\newline
e-mail:  \rm \texttt{glauco.valle@epfl.ch}
}
\date{}


\begin{document}

\begin{abstract} 
We show that for the voter model on $\{0,1\}^{\ZZ}$ corresponding to
a random walk with kernel $p(\cdot)$ and starting from unanimity to the right and opposing unanimity to
the left, a tight interface between $0$'s and $1$'s exists if $p(\cdot)$ has finite second moment but does
not if $p(\cdot)$ fails to have finite moment of order $\alpha$ for some $\alpha <2$.
\end{abstract}

\maketitle

%%%%%%%%%%%%%%%%%%%%%%%%%%%%%%%%%%%%%%%%%%%%%%%%%%%%%%%%%%%%%%%%%%%%%%%%%%%%
\section{Introduction}
%%%%%%%%%%%%%%%%%%%%%%%%%%%%%%%%%%%%%%%%%%%%%%%%%%%%%%%%%%%%%%%%%%%%%%%%%%%%
\setcounter{equation}{0}
\label{sec:introd}

The one-dimensional voter model is a spin system with configuration space $\Omega=\{0,1\}^{\ZZ}$ and flip rates
$$
c(x, \eta)= \sum_{y} p(y-x) \textbf{1}_{\eta(x)\neq \eta(y)},
$$
where $p(\cdot)$ is a transition probability kernel on $\ZZ$.  This kernel is taken to be irreducible, i.e., $\{x:p(x)+p(-x)>0\}$ generates $\ZZ$. The standard interpretation is that the values $0$ and $1$ represent two opinions and each site $x$
of $\ZZ$ represents an individual. For a configuration $\eta \in \Omega $, the value
$\eta(x)$ gives the ``opinion" of the individual $x$. Then the evolution of the system can be described in the following way: each individual, independently of the others, decides at times occurring as a Poisson process to adopt the opinion of a neighbor whose position is chosen according to the kernel $p(\cdot)$.

The voter model was introduced independently by \cite{HL} and \cite{Cli} where it was shown that 
non trivial extremal equilibrium measures exist (i.e. distinct from  the dirac measures $\delta_{\underline
0}, \delta_{\underline 1} $) if and only if the random walk corresponding to the symmetrized
probability $\tilde p(x) = (p(x)+p(-x))/2$ is transient.  In this latter case there exists an
ergodic equilibrium for every density between $0$ and $1$ and these formed the extremal
equilibria. On the other hand, for instance, in one-dimension if $\sum_x |x|p(x) \ < \ \infty $, then this symmetrized random walk is recurrent and there are only the trivial (unanimous) measures. We refer to \cite{L} and \cite{D1} for the statement and proof of these results and also for results on existence and uniqueness of the process and all relevant background material.

Let $\eta_{0}$ be the Heavyside configuration on $\Omega$, i.e., the configuration:
$$
\eta_{0}(z) =
\left\{ \begin{array}{ll}
1 , & \textrm{if } z \le 0 \\
0 , & \textrm{if } z \ge 1 \, ,
\end{array}\right.
$$
and consider the voter model $\left(\eta_t : t\geq 0 \right)$ starting at $\eta_{0}$.
For each time $t>0$, let 
$$
r_t = \sup\{x: \eta_t(x) =1 \} \,  \quad \hbox{ and } \quad
l_t = \inf\{ x: \eta_t(x) =0 \},
$$
which, when finite, are respectively the positions of the rightmost one and the leftmost zero. We call the voter model
configuration between the coordinates $l_t$ and $r_t$ the {\it voter model interface},
and $r_t-l_t+1$ is the {\it interface size}. From now on we suppose that the
probability kernel $p\,(\cdot)$ has finite absolute first moment which turns out to be a necessary and sufficient condition for the interfaces to be almost surely finite for all $t\ge 0$ and thus well defined. To see sufficiency, we observe that the rate at which the interface size increases is bounded above by
$$
\sum_{x<0<y} \{ p\,(y-x) + p\,(x-y) \} = \sum_{z\in \ZZ } |z| p\,(z) \, .
$$

Our starting point is Cox and Durrett's paper \cite{CD} where among other results it was shown
that if the voter model corresponds to a kernel $p(\cdot)$ with $\sum_x |x|^3 p(x) < \infty $, then
the voter model interface is tight, i.e., the random variables $(r_t - l_t)_{t \ge 0}$ are tight.  We use much of the overall approach of \cite{CD}, though much elegance is sacrificed, in proving

\medskip

\begin{thm} \label{thm1} {
For the one-dimensional voter model $\left(\eta_t : t\geq 0 \right)$ corresponding to an
irreducible kernel $p(\cdot)$ and starting from $\eta_0$, the random variables $(r_t -
l_t)_{t \geq 0}$ are tight if $\sum_x x^2 p(x) < \ \infty $.}
\end{thm}

\medskip

This result is almost optimal as our next result shows

\medskip

\begin{thm} \label{thm2} {
For the one-dimensional voter model $\left(\eta_t : t\geq 0 \right)$ corresponding to an
irreducible kernel $p(\cdot)$ starting from $\eta_0$, the random variables $\left(r_t -
l_t : t \geq 0\right)$ are not tight if $\sum_x |x|^c p(x) \ = \ \infty $ for some $c < 2$.}
\end{thm}

\medskip

\begin{remark} It was noted in \cite{CD} that the question of the tightness interface can be treated as
the question of the positive recurrence of the irreducible
Markov chain obtained from the voter model $\left(\eta_t : t\geq 0 \right)$ viewed from its leftmost zero, this is a Markov chain whose state space is the subset of $\{0,1\}^{\ZZ_+}$ having only finitely many ones. What is surprising is that, the trivial nearest neighbour case excepted, if this chain is positive recurrent, then, under the equilibrium
distribution, the position of the rightmost one is of infinite mean. Sharp estimates for the tail distribution are given in \cite{BMSV}. It was also observed in \cite{CD} that if the interface is tight then the finite dimensional distributions of the diffusively rescaled position of the rightmost one and leftmost zero converge to those of a Brownian motion with speed $\sigma=\sum x^2 p(x)$. The question of the convergence of diffusively rescaled position of the rightmost one and leftmost zero as a process on the space of continuous time right continous with left limit trajectories to a $\sigma$-speed Brownian motion was also considered in \cite{BMSV}.
\end{remark}

\medskip

We now give the Harris construction for the process and introduce its dual which
will play a crucial role in our analysis. Let $\{ \mathcal{N}^{x,y} \}_{x,y \in \ZZ}$ be independent
Poisson point processes with intensity $p(y-x)$ for each $x,y \in \ZZ$.
From an initial configuration $\eta_0$ in $\Omega$, we set at
time $t \in \mathcal{N}^{x,y}$:
$$
\eta_t(z) =
\left\{ \begin{array}{ll}
\eta_{t-}(z) , & \textrm{if } z \neq x \\
\eta_{t-}(y) , & \textrm{if } z = x \, .
\end{array}\right.
$$
This is unambiguous as a.s. no time $t$ can be in two distinct independent Poisson processes and each Poisson
process has no double points.  It is to be shown that this indeed defines a process and that there
is no problem of infinite regress. We refer the reader to \cite{D1} for details.

From the same Poisson point processes, we construct the system of \emph{coalescing random walks}
as follows. We can think of the Poisson points in ${\mathcal{N}}^{x,y}$ as marks at site $x$ occurring
at the Poisson times. For each space-time point $(x,t)$ we start a random walk $X^{x,t}$ evolving
backward in time such that whenever the walk hits a mark in ${\mathcal{N}}^{u,v}$ 
(i.e., for $s \in (0,t)$, $(t-s)\in{\mathcal{N}}^{u,v}$ and $u=X^{x,t}_s$), it jumps from site $u$
to site $v$. When two such random walks meet, which occurs because one walk
jumps on top of the other walk, they coalesce into a single random walk starting from the space-time
point where they first met. We define by $(\xi_s:s\ge 0)$ the Markov process which describes the positions
of the coalescing particles at time $s$. If $\xi_s$ starts at time $t$ with one particle from every
site of $A$ for some $A\subset \ZZ$, then we use the notation
$$
\xi^t_s(A) := \{ X^{x,t}_s : x \in A \} \, ,
$$
where the superscript is the time in the voter model when the walks first started, and the subscript
is the time for the coalescing random walks. It is well known that $\xi_t$ is the dual
process of $\eta_t$ (see Liggett's book \cite{L}), indeed the random walk $X_{s}^{x,t}$ traces backward in time the origin of the value (or opinion) at site $x$ at time $t$ in the sense that $\eta_t(x)=\eta_0(X_{t}^{x,t})$ and thus
\begin{equation}
\label{eq:dualrel}
\{ \eta_t(z) \equiv 1 \textrm{ on } A \} = \{ \eta_0(z) \equiv 1 \textrm{ on } \xi^t_t(A) \} \,
\end{equation}
for all $A\subset \ZZ$.

Since a finite family of dual coalescing random walks behaves as independent random walks until the first time two of them meet each other, it will be necessary all along the paper to estimate events related to the evolution of a system of independent random walks starting on each integer site. So we fix a notation, we will denote by $(Y^x:x\in \ZZ)$ a family of independent one-dimensional homogeneous translation invariant continuous time random walks jumping at rate 1 according to the transition probability kernel $p(\cdot)$ where the parameter $x$ indicates the starting point of the walk. 

\medskip

We are going to suppose from this point without loss of generality that the probability kernel $p(\cdot)$ has zero mean, i.e., $\sum_x x p(x)= 0$. We do this in part due to the fact that duality arguments reduce most of our proofs to the analysis of the behaviour of the difference between two independent random walks, which is necessarily a symmetric random walk. Beyond that the adjustments to be made in passing from the mean zero case to an arbitrary mean really just involve straightforward translations and changes of notation.

\medskip

The Proof of Theorem \ref{thm2} is entirely contained in Section \ref{sec:theorem2}. The other sections are devoted to the proof of Theorem \ref{thm1} and we end this section by discussing its proof. 

\medskip

The proof for tightness of the voter model interface under absolute third moment assumption in \cite{CD} was based on duality, based on the fact that  to understand the behaviour of the interface it is enough to understand the evolution of the dual coalescing random walks. Recall that tightness can be stated as
\begin{equation}
\label{eq:tight}
\lim_{M \ra +\oo} \sup_{t>0} \mathrm{P} ( r_t-l_t >M ) = 0 \, .
\end{equation}
We are considering the voter model $\eta_{\cdot}$ starting from the heavyside configuration $\eta_0$. Then we regard the evolution of the dual coalescing random walks system. For the event $\{r_t-l_t >M\}$ to happen there must exist two sites $x<x+M \le y$ such that $\eta_t(x)=0$ and $\eta_t(y)=1$ at time $t$. Thus the two coalescing random walks $X^{x,t}$ and $X^{y,t}$ must end up with $X_t^{x,t} \geq 0 > X_t^{y,t}$ which also means that they have not hit each other and as a consequense they behave as independent random walks during the dual time interval $[0,t]$. The basic idea in \cite{CD} was then to count these pairs of dual coalescing random walks and estimate the probability that for any of these pairs a crossing of the walks as formaly described above occurs, i.e.,
\begin{eqnarray*}
\mathrm{P} \left( \eta_t(x)=0,\eta_t(y)=1 \right) &=& \mathrm{P} \left(
X_{t}^{x,t}\geq 0> X_{t}^{y,t}\right) \\
&=& \mathrm{P} \left( \left\{ Y^{x}(t)\geq 0 > Y^{y}(t)\,\,\right\}\cap
\left\{ \forall \, s \in [0,t]\, , Y^{x}(s) \neq Y^{y}(s)  \right\}\right).
\end{eqnarray*}

To establish tightness it is only necessary to consider large values of $t$ at the probability in (\ref{eq:tight}). Then in Cox and Durrett's proof the voter model was considered at a time $t >> K$ with $K$ being a large fixed constant, from where the proof was divided in two parts. The first part was related to the evolution of the dual random walks that have started at time $t$ during the (dual) time interval $[K,t]$.  Let us first remark that the density of such walks over the integers will be small for $K$ large due to coalescence, indeed the density of coalescing random walks tends to zero as time tends to infinity. Let $\Gamma$ be the event that there exist two dual random walks that started at time $t$ and are neighbors at dual time $K$, i.e., with no other dual random walk between them at such time, which
\begin{enumerate}
\item[(i)] do not coalesce but cross each other over the next $t-K$ units of time and
\item[(ii)] at time $0$, the leftmost walk's position is in $(- \infty , 0)$ and the rightmost walk's position is in $[0, \infty )$.
\end{enumerate}
It was shown that (uniformly over $t > K$) the probability that $\Gamma$ happens is small if $K$ is large.

The second part consisted of showing that, if event $\Gamma$ occurs then the size of the interface would
have its expectation bounded by a constant depending on $K$ uniformly over $t > K$.  This step only
requires finite second moment and so we have no need to tighten it here.

\medskip

We prove Theorem \ref{thm1} extending the first part in the approach of \cite{CD} using only second moment assumptions. 
To explain the idea behind the proof, we start by introducing some notation. We define $\tau_A^S=\inf \{t > 0\, :\; S(t) \in A\}$ and $\widehat{\tau}_A^S=\inf \{t>\sigma_1\,:\; S(t) \in A \}$ for a continuous time stochastic process $(S(t): \ t \ \geq 0)$ and $A \subset \ZZ$, where $\sigma_1$ is the time of the first jump. If $S$ is given and it is unambiguous the suffix $S$ will be dropped, and if $A=\{x\}$, a single point, we use the notation $\tau_x$ and $\widehat{\tau}_x$ instead of $\tau_{\{x\}}$ and $\widehat{\tau}_{\{x\}}$ respectively. 

The idea is to divide up the sets of neighbouring dual random walks according to whether the position of the random walk nearest the origin at a large fixed dual time $K$ is of magnitude less then $m \sqrt t$ or whether it is larger than this. For pairs of magnitude less than $m \sqrt t$ we consider the event that the two random walks reverse order but do not coalesce. Ignoring the ultimate positions of the random walks which should only cost us a factor (depending on $M$) by the invariance principle. Here the main work is to estimate the probability of the event $V^k(t)= \{X^0(t)> X^k(t);\, \tau_0 >t\}$, where $\tau_0 =\tau_{0}^{X^k-X^0}$. We show that

\medskip

\begin{prop} 
\label{lemm1}
Let $V^k(t)$ be the event defined as above, then
$$ 
\lim_{l\rightarrow \infty} \limsup_{t\rightarrow \infty} \sup_{k>l} \frac{\mathrm{P}(V^k(t))}{k/\sqrt{t}}=0.
$$
\end{prop}

\medskip

For the pairs farther away from the origin we consider instead the event that the two random walks do not coalesce and that at least one of them crosses the origin.  Here the following result is the most relevant

\medskip

\begin{prop} 
\label{lemm2}
For a negative integer $k$, let $U_k(t)$ be the event $\{Y^{0}(s) \neq Y^{k}(s), \textrm{ for every, } 0 \leq s \leq t \}$. Then, there exists finite constant $c_0 > 0$, depending only on $p(\cdot)$,  such that for all negative $k$ and positive $m$
and $t>0$
$$
\mathrm{P}\left( U_k(t) \cap \left\{ Y^0(t)\geq m \sqrt{t} \;or\; Y^k(t)\geq m \sqrt{t}\right\}  \right)\leq
c_0 \mathrm{P}(U_k(t)) \mathrm{P}\left(|Y^0(2t)|)\geq m \sqrt{t}\right).
$$
\end{prop}

\medskip

We still fix here some notation that will be used all along the paper. For two sequences of real numbers $(A_n)_{n \ge 1}$ and $(B_n)_{n \ge 1}$, if there exists constants $b_1>0, \,\,b_2>0$ not depending on $n$ and $n_0>0$ such that for every $n>n_0$ we have $b_1 B_n \leq A_n \leq b_2 B_n $, we denote $A_n \approx B_n $. Given two collections of sequences indexed by $k$, $(A_{n}^{k})_{n \ge 1}$ $(B_{n}^{k})_{n \ge 1}$ we say that $A_{n}^{k}\approx B_{n}^{k}$ uniformly in $k$ if for
each fixed $k$, $A_{n}^{k}\approx B_{n}^{k}$ and the relevant $n_0,\,b_1,\, b_2$ can be chosen
independent of $k$. Unless otherwise stated, constants $c_i$ will all be positive, finite constants which may depend only on $p(\cdot)$. Moreover when dealing with some stochastic process on the integers, $\mathrm{P}^x$ will denote the probability when the process begins at the site $x$, and for a joint process $(X_1(\cdot), X_2(\cdot))$, $\mathrm{P}^{x,y}$ denotes the probability when the first process begins at the site $x$, and the second at the site $y$.

\smallskip

As we already mentioned the proof of Theorem \ref{thm2} is the content of Section \ref{sec:theorem2}. Sections \ref{sec:lemma1}, \ref{sec:lemma2} and \ref{sec:hybzones} are respectively devoted to the proofs of Proposition \ref{lemm1},
Proposition \ref{lemm2} and Theorem \ref{thm1}. Moreover an appendix have been included at the end of the paper containing results and properties of random walks needed in other sections.

%%%%%%%%%%%%%%%%%%%%%%%%%%%%%%%%%%%%%%%%%%%%%%%%%%%%%%%%%%%%%%%%%%%%%%%%%%%%%%
\section{Proof of Proposition \ref{lemm1}}
%%%%%%%%%%%%%%%%%%%%%%%%%%%%%%%%%%%%%%%%%%%%%%%%%%%%%%%%%%%%%%%%%%%%%%%%%%%%%%
\label{sec:lemma1}
\setcounter{equation}{0}

Recall $Y^0(\cdot)$ and $Y^k(\cdot)$ are two independent random walks starting at $0$ and $k>0$ respectively that have the same transition probability kernel $p(\cdot)$ which has finite second moment. The Proposition \ref{lemm1} is clearly a result on 
the symmetric continuous time random walk $(Y^k(t)-Y^0(t)\,:\, t \geq 0)$ whose transition kernel is $p(x)+p(-x)$, $x\in \ZZ$.
Moreover the act of rescaling the jump rates of $(Y^k(t)-Y^0(t)\,:\, t \geq 0)$ by a fixed constant does not change the result in Proposition \ref{lemm1} nor its proof. Then we can simply recast the statement of the Lemma stating that for a symmetric homogeneous translation invariant continuous time random walk $(Z(t)\,:\, t \geq 0)$ whose transition kernel has finite second moment
\begin{equation}
\label{eq:lemma1}
\lim_{l\rightarrow \infty} \limsup_{t\rightarrow \infty} \sup_{k>l} \frac{\mathrm{P}^{k} ( Z(t)<0, \tau^Z_0 >t)}{k/\sqrt{t}}=0.
\end{equation}

From now on we fix the symmetric homogeneous translation invariant continuous time random walk $(Z(t)\,:\, t \geq 0)$ associated to a transition probability kernel $\tilde{p}$ which has finite second moment. During this section $\tau_A$ refers to the hitting time for the random walk $Z$.

\medskip
 
The intuitive ideas behind the limiting behavior described in (\ref{eq:lemma1}) and its proof come from  
Lemma \ref{lem1} and the invariance principle{\footnote{The invariance principle is used many times during the paper. Although our applications of it are standard facts for mean zero finite variance random walks, we leave a reference here which is the text of Kruglov \cite{K} that gives a survey on the subject.}} and are based in the following argument: If Z starts at $k$ for some large $k>0$ then by Lemma \ref{lem1} we expect the probability of $\{Z(t)<0\}$ conditioned on $\{Z(s) \neq 0, 0<s<t\}$ to be of order smaller than $\frac{|Z(\tau_{(-\oo,0)})|}{\sqrt{t-\tau_{(-\oo,0)}}}$. By the invariance principle, $|Z(\tau_{(-\oo,0)})|$
should be small compared to $k$ thus $\frac{|Z(\tau_{(-\oo,0)})|}{\sqrt{t}}<< \frac{k}{\sqrt{t}}$. Two issues now present themselves.  First we must in fact deal with $E^k[Z(\tau_{(-\oo,0)})]$ and secondly since the conditional probability in question is actually of the order $\frac{Z(\tau_{(-\oo,0)})}{\sqrt{t- \tau_{(-\oo,0)} }}$, we must handle the case where $t- \tau_{(-\oo,0)} << t$.

\medskip

We show (\ref{eq:lemma1}) by considering separately the cases $\tau_{(-\infty, 0]} < t/2$ and $\tau_{(-\infty, 0]} \ge t/2$. Then, denoting by $\tilde{V}^k (t)$ the event $\{Z(0)=k, \, Z(t)<0 , \, \tau_0>t\}$, (\ref{eq:lemma1}) follows obviously from the next two results:

\medskip

\begin{lemma} \label{Vlemma}
Let $\tilde{V}^k(t)$ be the event defined as above, then
$$ 
\lim_{l\rightarrow \infty} \limsup_{t\rightarrow \infty} \sup_{k>l} \frac{\mathrm{P} \left(\tilde{V}^k (t) \cap \{ \tau_{(-\infty, 0]} < t/2\}\right)}{k/\sqrt{t}} = 0. 
$$
\end{lemma}
\no and
\begin{lemma} \label{Blemma}
Let $\tilde{V}^k(t)$ be the event defined as above, then
$$ 
\lim_{l\rightarrow \infty} \limsup_{t\rightarrow \infty} \sup_{k>l} \frac{\mathrm{P} \left(\tilde{V}^k (t) \cap \{ \tau_{(-\infty, 0]} > t/2\}\right)}{k/\sqrt{t}} = 0. 
$$
\end{lemma}

\bigskip

Let us start by fixing some notation. Let $I_{k,r}$ be the interval $(0,k2^r)$, then put 
$$
A(k,r):=\{ \tau_{I_{k,r}^c}< \tau_{(-\infty,0]}\,;\,\tau_{(-\infty,0]} \leq \tau_{I_{k,r+1}^c}\} \, ,
$$
and
$$ 
H_{k,r}(l) := \int_{0}^{\tau_{I_{k,r}^c}} I_{Z(s)=l} ds,
$$
with $H_{k,r}(l)$ so defined, for every $k>0$, $r\ge 0$ and $x\in \ZZ$ the expectation $E^x [H_{k,r}(l)]$ is the green function for the process $Z$ starting at $x$ and killed on leaving $I_{k,r}$.

\medskip
Before we give the proofs of Lemma \ref{Vlemma} and of Lemma \ref{Blemma} we will establish several estimates related to the events $A(k,r)$ and the green functions $E^x[H_{k,r}(l)]$. 

\medskip

\begin{lemma} 
\label{lem_04}
There exists finite constant $c_1 > 0$ so that for all positive integers $k$ and $r$, for every $(x,l) \in I_{k,r} \times I_{k,r}$,
$$
E^{x}\left(H_{k,r}(l)\right)=\frac{\mathrm{P}^x(\tau_{l}<\tau_{I_{k,r}^c})}{\mathrm{P}^l(\tau_{I_{k,r}^c}< \widehat{\tau}_l)}
\leq c_1 (x \wedge l).
$$
\end{lemma}

\medskip

\no \textbf{Proof:} Since the jump rate of the random walk $Z$ is $1$ the green function $E^x[H_{k,r}(l)]$ is indeed the expectation of the number of visits of the walk to site $l$ before it leaves the interval $I_{k,r}$. Denote $n_{k,r}(l)$ the random number of such visits. Therefore,
$$
E^{x}\left[ H_{k,r}(l) \right] = E^x [n_{k,r}(l)] = \sum_{n=0}^{\infty} n \mathrm{P}^{x}(n_{k,r}(l)=n) \, .
$$
By the strong Markov property
$$
\mathrm{P}^{x}(n_{k,r}(l)=n) = \mathrm{P}^x(\tau_l < \tau_{I_{k,r}^c}) \mathrm{P}^l(\widehat{\tau}_{l}<\tau_{I_{k,r}^c})^{n-1} \mathrm{P}^l(\widehat{\tau}_l>\tau_{I_{k,r}^c}) \, ,
$$
Hence
$$
E^{x}\left[ H_{k,r}(l) \right] = \sum_{n=0}^{\infty} n \mathrm{P}^x(\tau_l < \tau_{I_{k,r}^c}) \mathrm{P}^l(\widehat{\tau}_{l}<\tau_{I_{k,r}^c})^{n-1} \mathrm{P}^l(\widehat{\tau}_l>\tau_{I_{k,r}^c})= \frac{\mathrm{P}^x(\tau_l < \tau_{I_{k,r}^c})}{\mathrm{P}^l(\tau_{I_{k,r}^c}<\widehat{\tau}_l)},
$$
which gives the equality in the statement, moreover the rightmost hand term is clearly bounded by
$$
\frac{\mathrm{P}^x(\tau_l < \tau_0)}{\mathrm{P}^l(\tau_0<\widehat{\tau}_l)}.
$$
Now apply Lemma \ref{lem1}.\textbf{d} to the numerator and Lemma \ref{lem1}.\textbf{a} to the denominator in the fraction to conclude. $\square$

\medskip

\begin{lemma} 
\label{lemm03}
There exists a finite constant $c_2 > 0$, such that for all positive integers $k$ and $r$, $P^k(A(k,r))\leq c_2/2^r.$
\end{lemma}

\medskip

\no \textbf{Proof:} Since $A(k,r) \subset  \{\tau_{[k2^r, \infty)}<\tau_{(-\infty,0]}\}$, it suffices to prove the result obtained by replacing $A(k,r)$ by $\{\tau_{[k2^r, \infty)}<\tau_{(-\infty,0]}\}$ in the statement. Now note that 
\begin{eqnarray}
\label{eq:lemm03}
\lefteqn{ \mathrm{P}^k(\tau_{[k2^r, \infty)}<\tau_{(-\infty,0]})\mathrm{P}^{k2^r}(\tau_{(-\infty,0]}>(k2^r)^2) \le } \nn \\
& & \le \sum_{i=k2^r}^{\infty} \mathrm{P}^k(\tau_{[k2^r, \infty)}=\tau_i<\tau_{(-\infty,0]}) \mathrm{P}^i(\tau_{(-\infty,0]}>(k2^r)^2) \nn  \\
& & \le \mathrm{P}^k(\tau_{(-\infty,0]}\geq (k2^r)^2) \le \mathrm{P}^k(\tau_{0}\geq (k2^r)^2) ,
\end{eqnarray}
where the second inequlity comes from the strong Markov property. By the invariance principle, the probability $\mathrm{P}^{M}(\tau_{(-\infty,0]}>M^2)$ is bounded away from the origin by a strictly positive constant uniformly on $M \in \ZZ_+$. Therefore, by (\ref{eq:lemm03}), the previous observation and Lemma \ref{lem1}.\textbf{c}
\begin{equation}
\label{eq:rexit}
\mathrm{P}^k(\tau_{[k2^r, \infty)}<\tau_{(-\infty,0]}) \leq \frac{\mathrm{P}^k(\tau_{0} \ge (k2^r)^2)}{\mathrm{P}^{k2^r}(\tau_{(-\infty,0]}>(k2^r)^2)}
\le  \frac{ c_2 k}{ k2^r} = \frac{c_2}{2^r},
\end{equation}
for some constant $c_2>0$. $\Box$ 

\bigskip

This bound is in fact an equivalence

\medskip

\begin{lemma} 
\label{lemm13}
Uniformly in $k$, the events $A(k,r)$ defined above satisfy, $\mathrm{P}^k(A(k,r))\approx 1/2^r.$
\end{lemma}

\medskip

\no \textbf{Proof:} We have already show that $\mathrm{P}^k(A(k,r))\leq c_2/2^r$ in Lemma \ref{lemm03} for some finite $c_2$, so
to prove this lemma it suffices to show that there exists $c_3>0$ such that for all positive integers $k$ and $r$ we have that $\mathrm{P}^k(A(k,r)) \geq c_3/2^r$. We make three claims from where the result follows:

\medskip

\no \textbf{Claim 1:} For all positive integers $k$ and $r$, we have that 
$$\mathrm{P}^k(A(k,r))\geq 4^{-1} P^k\left(Z({\tau_{I_{k,r}^c}}) \in [k2^r, 3k2^r/2]\right).$$

\medskip

\no \textbf{Claim 2:} For all positive integers $k$ and $r$, we have that  
$$\mathrm{P}^k\left(Z({\tau_{I_{k,r}^c}}) \geq k2^r\right)\geq 1/2^r \, .$$

\medskip

\no \textbf{Claim 3:} $2^r \mathrm{P}^k\left(Z({\tau_{I_{k,r}^c}})> 3k2^r/2\right)\rightarrow 0$ as $k+r$ tends to infinity.

\medskip

To conclude we apply the claims to obtain that
\begin{eqnarray*}
  \mathrm{P}^k(A(k,r)) & \ge & c \, 4^{-1} \mathrm{P}^k \left(Z({\tau_{I_{k,r}^c}}) \in [k2^r, 3k2^r/2] \right) \\
    &\geq& c \, 4^{-1} \left( \mathrm{P}^k\left(Z({\tau_{I_{k,r}^c}}) \geq k2^r\right)-\mathrm{P}^k\left(Z({\tau_{I_{k,r}^c}})> 3k2^r/2\right) \right) \\
    &\geq& c_2/2^r
\end{eqnarray*}
for all but a finite number of pairs $(k,r)$ and the result now follows by the strict positivity of $\mathrm{P}^k(A(k,r))$.
So it remains to show these three claims.

\medskip

\no \textbf{Proof of Claim 1:} By the strong Markov property,
\begin{eqnarray*}
  \mathrm{P}^k(A(k,r)) &=& \sum_{l=k2^r}^{k2^{r+1}-1} \mathrm{P}^k\left(Z({\tau_{I_{k,r}^c}})=l\right) \mathrm{P}^l\left(Z({\tau_{I_{k,r+1}^c}})<0\right) \\
&\geq& \sum_{l=k2^r}^{3k2^{r}/2} \mathrm{P}^k\left(Z({\tau_{I_{k,r}^c}})=l\right) \mathrm{P}^l\left(Z({\tau_{I_{k,r+1}^c}})<0\right).
\end{eqnarray*}\
\newline
Since we have for all $k>1$, $r>0$ and $l\leq {3k2^r/2}$  that
$$
\mathrm{P}^l\left(Z({\tau_{I_{k,r+1}^c}})<0 \right) \geq \mathrm{P}^{3k2^r/2}\left(Z({\tau_{I_{k,r+1}^c}})<0 \right) \geq \frac{1}{4},  
$$
where the last inequality follows from the symmetry of $Z(\cdot)$, then
$$
\mathrm{P}^k(A(k,r)) \ge  4^{-1} \sum_{l=k2^r}^{3k2^{r+1}/2} \mathrm{P}^k\left(Z({\tau_{I_{k,r}^c}})=l\right)
\ge 4^{-1} \mathrm{P}^k\left(Z({\tau_{I_{k,r}^c}}) \in [k2^r, 3k2^r/2]\right).
$$

\medskip

\no \textbf{Proof of Claim 2:} As in the proof of Claim 1, we have that Claim 2 is a consequence of the symmetry of $Z(\cdot)$.

\medskip

\no \textbf{Proof of Claim 3:} The event that $\{Z(\tau_{I_{k,r}^c})> 3k2^r/2\}$ is equal to the disjoint union over sites $l$ from $1$ to $k2^r-1$, of the event that the exit of the random walk $Z_.$ from the inteval $I_{k,r}$ occurs in a jump from $l$ to $(3k2^r/2, \infty)$, that is to say a positive jump to the right from site $l$ of size greater than $ 3k 2^r/2-l$ occurs before the random walk has exited the interval $I_{k,r} $.  Since $Z_.$ jumps at rate $1$, then via a limit argument we see that for each $l \in (0,k2^r)$ this event is smaller or equal to
$$
E^k \left( \int_0 ^ \infty I_{Z(s) = l, \tau_{I_{k,r}^c} > s} \, ds \right) \sum_{x> 3k 2^r/2-l} \tilde{p}(x)
$$
and so
\begin{eqnarray*}
\lefteqn{ \!\!\!\!\!\!\!\!\!\!\!\!\!\!\!\!\!\!\!\!\!\!\!\!\!\!\!\!\!\!\!\!\!
2^r \mathrm{P}^k \left(Z(\tau_{I_{k,r}^c})> 3k2^r/2\right) \le  
2^r \sum_{l=1}^{k2^r-1} \mathrm{P}^k ( Z(\tau^c_{k,r}-)=l ) \left(\sum_{x> 3k 2^r/2-l}\tilde{p}(x)\right) } \, \\  
& & \le 2^r \sum_{l=1}^{k2^r-1} E^{k} \left[ H_{k,r}(l)\right] \left( \sum_{x> 3k 2^r/2-l} \tilde{p}(x)\right) \, .
\end{eqnarray*}
By Lemma \ref{lem_04} the right hand term is dominated by
$$
c_1 \, k 2^r \sum_{l=1}^{k2^r-1} \left( \sum_{x> 3k 2^r/2-l} \tilde{p}(x) \right) 
$$
which is bounded above by
$$
c_1 (k 2^r)^2 \left( \sum_{x> k 2^r/2} \tilde{p}(x) \right)
$$
that converges to zero as $k2^r$ tends to infinity, since the finite second moment assumption implies that $x^2 \sum_{|z|\geq x} \tilde{p}(z)$ tends to zero as $x$ tends to infinity. Thus the proof of the lemma is complete. $\Box$

\bigskip

\begin{lemma} \label{lemm23}
The following two limits hold:
$$
\lim_{k \ra +\oo} \sup_r \frac{ E^k\left[|Z({\tau_{(-\infty,0]}})|\; \big| \,  A(k,r)\,  \right]}{k 2^r} = 0 \quad \textrm{and} \quad 
\lim_{k \ra +\oo} \frac{ E^k\left[ |Z({\tau_{(-\infty,0]}})| \right]}{k} = 0 \, .
$$
\end{lemma}

\medskip  
  
\begin{remark}
The proof of Lemma \ref{lemm23} also allows us to conclude that $\sup_{k,r} E^k\left[|Z({\tau_{(-\infty,0]}})|\; \big| \,  A(k,r)\,  \right]$ and $\sup_{k} E^k\left[ |Z({\tau_{(-\infty,0]}})| \right]$ are finite if $\tilde{p}(\cdot)$ has finite absolute third moment.
\end{remark}

\medskip

\no \textbf{Proof:} We start by applying the strong Markov property and the symmetry of $\tilde{p}(\cdot)$ to decompose 
$E^k\left[|Z({\tau_{(-\infty,0]}})|,A(k,r)\right]$ as
$$
\sum_{y=1}^{k2^{r+1}-1} \sum_{m=0}^{\infty} \sum_{l=k2^r}^{k2^{r+1}-1} \mathrm{P}^k\left(Z({\tau_{I_{k,r}^c}})=l\right) E^{l}\left[ H_{k,r+1}(y)\right] m \tilde{p}(y+m) \, .
$$
Now apply Lemma \ref{lem_04} to see that this last expression is bounded above by 
\begin{eqnarray*}
\lefteqn{ \!\!\!\!\!\!\!\!\!\!\!\!\!\!\! c_1 \sum_{y=1}^{k2^{r+1}-1} \sum_{m=0}^{\infty} \left\{ \sum_{l=k2^r}^{k2^{r+1}-1} \mathrm{P}^k\left(Z({\tau_{I_{k,r}^c}})=l\right) \right\} y \, m\, \tilde{p}(y+m) \le } \\
& \le & c_1 \sum_{y=1}^{k2^{r+1}-1} \sum_{m=0}^{\infty} \mathrm{P}^k(\tau_{[k2^r, \infty)}<\tau_{(-\infty,0]}) \, y \, m\, \tilde{p}(y+m),
\end{eqnarray*}
which by (\ref{eq:rexit}) is dominated by
$$
\frac{c_1 \, c_2}{2^r} \sum_{y=1}^{k2^{r+1}-1} \sum_{m=0}^{\infty} y \, m \, \tilde{p}(y+m) ,
$$
Therefore
\begin{eqnarray}
\label{eq:ceA}
E^k\left[|Z({\tau_{(-\infty,0]}})|,A(k,r)\right] &
\le & \frac{c_1 \, c_2}{2^r} \sum_{z=1}^{\infty} \left(\sum_{y=0}^{z \wedge k2^{r+1}} y (z-y)\right) \tilde{p}(z) \nn \\
& \le & \frac{c_1 \, c_2}{2^r} \sum_{z=1}^{\infty} (z^3 \wedge k2^{r+1} z^2 \wedge k^2 2^{2(r+1)} z) \tilde{p}(z) \nn \\
& \le & \frac{c_1 \, c_2}{2^r} \sum_{z=1}^{\infty} (z^3 \wedge k^2 2^{2(r+1)} z) \tilde{p}(z) , 
\end{eqnarray}
for $c>0$ not depending on $k$ or $r$. So we have by Lemma \ref{lemm13} that there exists a constant $c>0$ such that
$$
\frac{E^k \left[ |Z({\tau_{(-\infty,0]}})| \big|  A(k,r) \, \right]}{k 2^r} =
\frac{E^k \left[ |Z({\tau_{(-\infty,0]}})| ,  A(k,r) \,  \right]}{k 2^r \, \mathrm{P}^k(A(k,r))} \leq c \, \sum_{z=0}^{\infty} \left(\frac{z^3}{k 2^r} \wedge z^2 \right) \tilde{p}(z),
$$ 
from where, since $\tilde{p}(\cdot)$ has finite second moment, we obtain that the first limit in the statement by the dominated convergence theorem. To show that that the second limit holds note that
$$
E^k[ | Z_{\tau_{(-\infty,0]}} | ] = \sum_{r\geq 0} E[|Z_{\tau_{(-\infty,0]}} |, A(k,r)]
$$
and then apply inequality (\ref{eq:ceA}) to bound the right hand term by
$$ 
\sum_{r=0}^{\infty} \sum_{z=1}^{\infty} \frac{c}{2^r} (z^3 \wedge k^2 2^{2(r+1)} z) \tilde{p}(z) \le c^{\prime} \sum_{z=1}^{\infty} (z^3 \wedge k z^2) \tilde{p}(z)
$$
for $c^{\prime}>0$ not depending on $k$ or $r$. To obtain the previous inequality simply note that
$$
2^{-r} \sum_{r:k2^r \le z} (k2^r)^2 z = k^2 z \sum_{r=1}^{l: k2^l \le z < k 2^{l+1}} 2^r \le k^2 z 2^{l+1} \le 2 k z^2 \, .
$$
Thus
$$
\frac{E^k[ | Z_{\tau_{(-\infty,0]}} | ]}{k} \le c \sum_{z=1}^{\infty} \left(\frac{z^3}{k} \wedge z^2 \right) \tilde{p}(z)
$$
which goes to zero as $k \ra +\oo$ again by the dominated convergence theorem.

The uniform bounds in the statement also follow directly from the previous computations.
$\Box$ 

\bigskip

Now we have all the results we need to show Lemma \ref{Vlemma}:

\bigskip

\no \textbf{Proof of Lemma \ref{Vlemma}:} We start by applying the strong Markov property to obtain the following decomposition
\begin{eqnarray*}
\lefteqn{
\mathrm{P}\left(\tilde{V}^k(t) \cap \{ \tau_{(-\infty, 0]}<t/2\}\right) \le } \\
    &\leq& \sum_{i=1}^\infty \mathrm{P}^k \left(\{ Z({\tau_{(-\infty,0]}})=-i\} \cap \{ \tau_{(-\infty, 0]}<t/2\}\right) \mathrm{P}^{\{-i\}}\left(\tau_{0}>\frac{t}{2}\right)
\end{eqnarray*}
By Lemma \ref{lem1}.\textbf{c}, the last expression is bounded above by 
$$
\gamma_3 \frac{\sum_{i=1}^\infty i \, \mathrm{P}^k \left(\{ Z({\tau_{(-\infty,0]}})=-i\} \right)}{\sqrt{t/2}} = \frac{\gamma_3}{\sqrt{t/2}} E\left[|Z({\tau_{(-\infty,0]}})|\right].
$$
Therefore we have that
$$
\frac{\mathrm{P}\left(\tilde{V}^k(t) \cap \{ \tau_{(-\infty, 0]}<t/2\}\right)}{\frac{ k}{\sqrt{t}}}\leq \sqrt{2}\, \frac{\gamma_3}{k} E^k \left[ | Z({\tau_{(-\infty,0]}}) | \right] \, ,
$$
and we finish the proof by applying Lemma \ref{lemm23}. $\Box$ 

\bigskip

It remains to deal with events $\tilde{V}^k(t) \cap \{ \tau_{(-\infty, 0]} > t/2\}$, where the random walk takes a ``long time" to reach the negative half line. But in order to show Lemma \ref{Blemma}, we will need one more result.

\medskip

\begin{lemma} 
\label{Blemma1} 
There exists $c_3>0$ such that
$$
\mathrm{P}^x\left( Z({(k2^r)^2})=y \; \big| \;\tau_{I_{k,r+1}^c}>(k2^r)^2 \right)\le c_3\, \frac{y \wedge
(k2^{r+1}-y)}{(k2^r)^2}
$$
for every $x,y \in I_{k,r+1}$.
\end{lemma}

\medskip

\no \textbf{Proof:} Let $w=k2^r$, for notational clarity. We also fix $x,y \in I_{k,r+1}$. Observe, as is clear from Lemma \ref{lem1} and Theorem 22.T1 in \cite{Spi}, that
\begin{equation}
\label{eq:p1}
\mathrm{P}^x(\tau_{I_{k,r+1}^c}>w^2)\approx \frac{x \wedge (2w-x)}{w}.
\end{equation}
We use the Markov property to obtain the following decomposition
\begin{eqnarray}
\label{eq:p2}
\lefteqn{ \mathrm{P}^x ( Z({w^2})=y , \tau_{I_{k,r+1}^c}>w^2)  = } \nn \\
& & = \sum_{l=1}^{k2^{r+1}-1} \mathrm{P}^x \left(\tau_{I_{k,r+1}^c}>\frac{w^2}{2}, Z\left(\frac{w^2}{2}\right)=l\right) \, \mathrm{P}^l \left(\tau_{I_{k,r+1}^c}>\frac{w^2}{2} ,Z\left(\frac{w^2}{2}\right)=y \right).
\end{eqnarray}
By symmetry the right hand side is equal to 
$$
 \sum_{l=1}^{k2^{r+1}-1} \mathrm{P}^x \left(\tau_{I_{k,r+1}^c}>\frac{w^2}{2}, Z\left(\frac{w^2}{2}\right)=l\right) \, \mathrm{P}^y \left(\tau_{I_{k,r+1}^c}>\frac{w^2}{2} ,Z\left(\frac{w^2}{2}\right)=l \right),
$$
and this sum is clearly bounded above by
\begin{equation}
\label{eq:p5}
\mathrm{P}^x \left(\tau_{I_{k,r+1}^c}>\frac{w^2}{2}\right) \,
\max_{l \in I_{k,r+1}} \mathrm{P}^y \left(\tau_{I_{k,r+1}^c}>\frac{w^2}{2}, Z\left(\frac{w^2}{2}\right)=l \right) \, .
\end{equation}
To deal with the second probability just above, we apply a decomposition similar to the one in (\ref{eq:p2}) to obtain that
\begin{eqnarray}
\label{eq:p4}
\lefteqn{ \mathrm{P}^y \left(\tau_{I_{k,r+1}^c}>\frac{w^2}{2}, Z\left(\frac{w^2}{2}\right)=l \right) } \nn \\
& \le & \sum_{v=1}^{k2^{r+1}-1} \mathrm{P}^y \left(\tau_{I_{k,r+1}^c}>\frac{w^2}{4}, Z\left(\frac{w^2}{4}\right)=v\right) \, \mathrm{P}^v \left(\tau_{I_{k,r+1}^c}>\frac{w^2}{4} ,Z\left(\frac{w^2}{4}\right)=l \right)
\end{eqnarray}
By the local central limit theorem, see for instance \cite{D1}, there exists a constant $c>0$ not depending on $l$, $v$ and $w$ such that
$$
\mathrm{P}^v \left(\tau_{I_{k,r+1}^c}>\frac{w^2}{4} ,Z\left(\frac{w^2}{4}\right)=l \right) \le \frac{c}{w}
$$
Therefore, from (\ref{eq:p4}), it follows that
$$
\max_{l \in I_{k,r+1}} \mathrm{P}^y \left(\tau_{I_{k,r+1}^c}>\frac{w^2}{2}, Z\left(\frac{w^2}{2}\right)=l \right) \le
\frac{c}{w} \mathrm{P}^y \left(\tau_{I_{k,r+1}^c}>\frac{w^2}{4}\right),
$$
and then, from (\ref{eq:p2}) and (\ref{eq:p5}), we arrive at
\begin{equation}
\label{eq:p3}
\mathrm{P}^x ( Z({w^2})=y , \tau_{I_{k,r+1}^c}>w^2) \le \frac{c}{w} \mathrm{P}^x \left(\tau_{I_{k,r+1}^c}>\frac{w^2}{2}\right)
\mathrm{P}^y \left(\tau_{I_{k,r+1}^c}>\frac{w^2}{4}\right) \, .
\end{equation}
Now, by Lemma \ref{lem1}.\textbf{c}, we have that for all $l\in I_{k,r+1}$, $b>0$
$$
\mathrm{P}^l \left(\tau_{I_{k,r+1}^c}>\frac{w^2}{b^2}\right)
  \leq  \mathrm{P}^l \left(\tau_{0}>\frac{w^2}{b^2}\right) \wedge \mathrm{P}^{l}\left(\tau_{2w}>\frac{w^2}{b^2}\right) \leq b \gamma_3 \frac{ l\wedge (2w-l)}{w} \, .
$$
From this last inequality and (\ref{eq:p3}) we conclude that there exists a constant $c^{\prime}>0$ not depending on $x$, $y$ or $w$ such that  
$$
\mathrm{P}^x ( Z({w^2})=y , \tau_{I_{k,r+1}^c}>w^2) \le c^{\prime} \frac{(x \wedge (2w-x)) ( y \wedge (2w-y) )}{w^3}
$$
Then by (\ref{eq:p1}) we obtain the inequality in the statement. $\Box$ 

\bigskip

Now we are ready to prove Lemma \ref{Blemma}.

\bigskip

\no \textbf{Proof of Lemma \ref{Blemma}:} For every $t>0$ and all positive integers $k$ and $r$, denote by $B(k,r,t)$ the event
$$
\left\{ A(k,r) \cap \tilde{V}^k(t) \cap \{ \tau_{(-\oo,0]} \ge t/2 \} \right\} \, .
$$
Then $\tilde{V}^k(t) \cap \{ \tau_{(-\infty, 0]} > t/2\}$ is the disjoint union of $B(k,r,t)$ over $r$ and the limit in the statement can be rewritten as
$$ 
\lim_{l\ra \oo} \limsup_{t\rightarrow \infty} \sup_{k>l} \frac{\sum_{r} \mathrm{P}(B(k,r,t))}{k/\sqrt{t} } = 0 \, .
$$

\medskip

Write
$$
\mathrm{P} (B(k,r,t)) = \int_{t/2}^{t} \mathrm{P}^k(A(k,r) \cap \tilde{V}^k(t) \cap \{\tau_{(-\infty,0]} \in ds\}) \, .
$$
The above integral is dominated by
\begin{equation}
\label{Blemma0-1}
\int_{\frac{t}{2}}^{t} \sum_{y=1}^{k2^{r+1}-1} \sum_{v \neq 0} \mathrm{P}^k\left(A(k,r), \tau_{(-\infty,0]} > s_-, Z({s_-})=y\right) \tilde{p}(v-y) \mathrm{P}^v(\tau_0>t-s)ds \, .
\end{equation}
To estimate the integrand in the previous expression we start by obtaining good bounds on
$$
\mathrm{P}^k\left(A(k,r) , \tau_{(-\infty,0]}>s , Z(s)=y \right)
$$
for $s > t/2$ and $y \in I_{k,r+1}$. We start by applying the strong Markov property to get, as an upper bound for the last probability, the following decomposition
$$
\sum_{v \in I_{k,r+1}^c} \mathrm{P}^k\left(A(k,r) , \tau_{(-\infty,0]}>\frac{t}{5} , Z\left(\frac{t}{5}\right)=v\right)
\mathrm{P}^v\left(\tau_{I_{k,r+1}}>s-\frac{t}{5} , Z\left(s-\frac{t}{5}\right)=y \right) \, .
$$
The second probability in the sum just above is bounded by $cy/(s-t/5)$ for some constant $c>0$, see the proof of Lemma \ref{Blemma1}, and then
\begin{equation}
\label{Blemma0-2}
\mathrm{P}^k\left(A(k,r) , \tau_{(-\infty,0]}>s , Z(s)=y \right) \le \frac{10 c}{3} \frac{y}{t} \mathrm{P}(A(k,r))
\le \frac{10 \, c \, c_2}{3} \frac{y}{2^r t} \, ,
\end{equation}
where the last inequality comes from Lemma \ref{lemm03}.

\bigskip

Now apply (\ref{Blemma0-2}) to the first probability in (\ref{Blemma0-1}) and Lemma \ref{lem1}.\textbf{c} to the second one, we have then 
\begin{eqnarray*}
\mathrm{P}(B(k,r,t))  &\leq&  \gamma_3 \frac{10 \, c \, c_2}{3} \int_{t/2}^{t}  \sum_{y=1}^{k2^{r+1}-1}\sum_{v > 0} \frac{y}{2^r t} \frac{v}{\sqrt{t-s}}\, \tilde{p}(v + y) ds,  \\
    &\leq& \gamma_3 \frac{10 \, c \, c_2}{3 \sqrt{2}} \frac{1}{2 ^r \sqrt{t}} \sum_{y=1}^{k2^{r+1}-1} \sum_{v > 0} y v \, \tilde{p}(v + y).
\end{eqnarray*}
Therefore, as in the proof of Lemma \ref{lemm23},
\begin{eqnarray*}
  \frac{\sum_r \mathrm{P}(B(k,r,t))}{k/\sqrt{t} } &\leq& c \sum_{z>0} \sum_{r \geq 0}\frac{1}{k2^r} \left( z^3 \wedge (k2^r)^2z \right) \tilde{p}(z) \\
    &\leq& c \sum_z  \left(\frac{|z|^3}{k} \wedge z^2 \right) \tilde{p}(z) \, .
\end{eqnarray*}
for some constant $c>0$. By the dominated convergence theorem the last term in the previous expression goes to zero as $k$ tends to infinity. $\square$

\bigskip

%%%%%%%%%%%%%%%%%%%%%%%%%%%%%%%%%%%%%%%%%%%%%%%%%%%%%%%%%%%%%%%%%%%%%%%%%%%%%%%%%%%%%%%%
\section{ Proof of Proposition \ref{lemm2}}\
%%%%%%%%%%%%%%%%%%%%%%%%%%%%%%%%%%%%%%%%%%%%%%%%%%%%%%%%%%%%%%%%%%%%%%%%%%%%%%%%%%%%%%%%
\label{sec:lemma2}
\setcounter{equation}{0}

We start this section by fixing some notation. Let $Y_1$ and $Y_2$ be two independent continuous time random walks with a mean zero transition probability kernel $p(\cdot)$ satisfying the conditions in Section \ref{sec:introd}. We denote by $Z$ the random walk $Y_1 - Y_2$, and in this section $\tau_A$ and $\widehat{\tau}_A$, for $A \subset \ZZ$, refer to the hitting times of the random walk $Z$. Note that the random walk $Z$ is of jump rate 2 instead of 1 and has transition kernel given by $p(x) + p(-x)$, $x\in \ZZ$. We will also be using that the superscript in $\mathrm{P}^{l,m}$ refers to the starting point of $(Y_1,Y_2)$ 
and in $\mathrm{P}^{l}$ it refers to the starting point of a single random walk which unless otherwise stated will be $Z$.

In such settings the inequality in the statement of Proposition \ref{lemm2} reads as
\begin{equation}
\label{eq:lemma2-1}
\mathrm{P}^{0,k}\left( \tau_0 > t \, , \, \max(Y_1(t), Y_2(t)) \geq m \sqrt{t} \right)\leq 
c_0 \, \mathrm{P}^{k} (\tau_0 > t) \, \mathrm{P}^0 \left(|Y_1(2t)|)\geq m \sqrt{t}\right)
\end{equation}
for all negative integers $k$, positive $m$ and $t>0$, where the event $\{\tau_0 > t\}$, for $(Y_1,Y_2)$ starting at $(0,k)$, plays the role of $U_k(t)$.

It is easily seen from the invariance principle that 
$$
\inf_{t > 0}P(Y_1(t) \geq 0 | Y_1(0) = 0) \quad \textrm{ and } \quad
\inf_{t > 0}P(Y_1(t) \leq 0 | Y_1(0) = 0) 
$$ are both strictly positive.  Hence by the strong Markov property we have the existence of $c > 0$ so that for all $t $,
$$
\mathrm{P}^0 \left(|Y_1(2t)|)\geq m \sqrt{t}\right) \geq c\mathrm{P}^0 \left( \sup_{s \leq 2t}|Y_1(s)|)\geq m \sqrt{t}\right).
$$
Hence it will suffice to show that
\begin{equation}
\label{eq:lemma2-11}
\mathrm{P}^{0,k}\left( \tau_0 > t \, , \, \max(Y_1(t), Y_2(t)) \geq m \sqrt{t} \right)\leq 
c_0 \, \mathrm{P}^{k} (\tau_0 > t) \, \mathrm{P}^0 \left( \sup _{s \leq 2t }|Y_1(s)|)\geq m \sqrt{t}\right).
\end{equation}

To show (\ref{eq:lemma2-11}), first observe that we only need to consider values of $t$ greater than $2 k^2$. Indeed, enlarging the constant $c_0$ if necessary, we have that (\ref{eq:lemma2-1}) also holds for $t\le 2k^2$ as a consequence of Lemma \ref{lem1}.\textbf{c} and the invariance principle which implies that
\begin{equation}
\label{eq:lemma2-2}
\mathrm{P}^{0,0} \left( \sup_{s\leq 2t} \max \left(|Y_1(s)|,\,|Y_2(s)|\right) \geq m\sqrt{t} \right)\leq 4
\mathrm{P}^0 \left( | Y_1(2t)| \geq  m\sqrt{t} \right).
\end{equation}
Hence we assume $t>2 k^2$ from now on.

\medskip

We will call an excursion, or excursion away from $0$, of $Z(\cdot)$ starting at time $w$ and having duration at least $t$ a random path of the type $((Z(s))_{w \le s \le \widehat{\tau}_0}: \ Z(w)=0, \ \tau_0 > t+w)$. Let $D_1$ be the event that an excursion (of $Z$) away from $0$ starts at some random time $s_1 \le t/2$ and that this excursion hits point $|k|$ before time $s_1 + k^2$ and thereafter lasts at least time $t$. Let $s_1\leq s_2 \leq s_1 +k^2$ be the random time that $Z(\cdot)$ first visits state $|k|$ after random time $s_1$. We also denote by $D_1^+$ the event that $D_1$  happens and $\max(Y_1(s_2), Y_2(s_2))=Y_1(s_2)>0$. By the invariance principle,
it can easily be seen, there exists a constant $c>0$ such that $\mathrm{P}^{0,0}(D_1^+) \ge c^{-1} P^0(D_1)$.

Given the event $D_1^+$, the process $\left( Y_1(s+s_2)-Y_1(s_2), Y_2(s+s_2)-Y_1(s_2):\, 0\leq s\le t \right)$ are equal in law to $\left(Y_1(s), Y_2(s):\, 0\leq s\leq t \right)$ starting at $(0,k)$ and conditioned on the event $\{\tau_0 > t\}$. Thus we have
\begin{eqnarray*}
\mathrm{P}^{0,k}\left( \max(Y_1(t), Y_2(t)) \geq m \sqrt{t} \, \big| \, \tau_0 > t \right)
& \le & \mathrm{P}^{0,0}\left( \sup_{0\le s \le 2t} \max_{i=1,2} Y_i(s) \geq m \sqrt{t} \, \big| \, D_1^+ \right) \\
& \le & \frac{\mathrm{P}^{0,0}\left( \sup_{0\le s \le 2t} \max_{i=1,2} Y_i(s) \geq m \sqrt{t} \, , \, D_1^+ \right)}{\mathrm{P^{0,0}}(D_1^+)} \\
& \le & c \, \frac{\mathrm{P}^{0,0}\left( \sup_{0\le s \le 2t} \max_{i=1,2} Y_i(s) \geq m \sqrt{t} \right)}{\mathrm{P^{0}}(D_1)} \,
\end{eqnarray*}
where the supremum over $0 \le s \le 2t $ is appropriate since $s_2 + t < 2t$. Now, if we can show that there exists a constant $h>0$ such that $\mathrm{P}(D_1)\geq h>0$ uniformly in $k$ and $t$ with $t > 2 k^2$, then (\ref{eq:lemma2-11}) follows from this last inequality by (\ref{eq:lemma2-2}).

\medskip

We end this section establishing the previous uniform lower bound and therefore completing the proof of (\ref{eq:lemma2-1}) and Proposition \ref{lemm2}. Let $E_1$ be the event that an excursion away from $0$ of $Z(\cdot)$ of duration at least $3t/2$ starts at some random time smaller than $t/2$. Note that $D_1 \cap E_1$ differs from $D_1$ in the fact that excursions in the former have duration greater than a fixed time lenght of $3t/2$, while the excursions from paths in $D_1$ have duration of at least the random time $s_2 + t \in (t,3t/2)$. Write
$$
\mathrm{P}^0 (D_1) \ge \mathrm{P}^0 (D_1 \cap E_1) \ge \mathrm{P}^0 (D_1 \big| E_1) \mathrm{P}^0 (E_1) \, ,
$$
and the desired result follows from the next two lemmas:

\medskip

\begin{lemma} 
\label{lemm04}
Let $E_1$ be the event defined just above. Then there exists $h_1>0$ such that for all $t>2$ we have $\mathrm{P}^0(E_1) \ge h_1$.
\end{lemma}

\medskip

\begin{lemma} \label{cor2}
Let $E_1$ be the event defined just above. Then, conditioned on $E_1$, the probability that the excursion reaches
$k$ before time $k^2$ and thereafter lasts time $t$ is bounded below by a constant $h_2>0$ uniformly over $k$ and $t> 2 k^2$.
\end{lemma}

\bigskip

We prove these two Lemmas in their order:

\bigskip

\no \textbf{Proof of Lemma \ref{lemm04}:} We are going to prove that the result holds for all $t$ sufficiently large. It is clear that by enlarging the constant if necessary we also have it for all $t>2$. Now let $(B(t): t\ge 0)$ be a standard Brownian motion. By the invariance principle we have that
$$
\mathrm{P}^0(E_1) \ge \mathrm{P}(\forall s \in [1/2, 2]: B(s) \neq 0) - \epsilon_a
$$
for all $t>a$, where $\epsilon_a \rightarrow 0$ as $a\rightarrow \infty$. We have
\begin{eqnarray*}
\mathrm{P}^0(\forall s \in [1/2, 3/2]: B_s \neq 0)	&=& 2 \int_0^{\infty} \mathrm{P}^0(B_{1/2} \in dx) \mathrm{P}^x(\tau_0^B >1)\\
&\geq& 2 \int_1^{\infty} \frac{e^{-x^2}}{\sqrt{\pi }} dx \mathrm{P}^1(\tau_0^B > 1) 
\end{eqnarray*}
and the proof is complete.$\Box$

\bigskip

\no \textbf{Proof of Lemma \ref{cor2}:} By the strong Markov property, the conditional probability described in the statement of the Lemma is equal to $\mathrm{P}^0 ( \tau_k < k^2 \, | \, \widehat{\tau}_0 > 3t/2)$. This probability is equal to
\begin{eqnarray*}
 & \ &
\frac{\mathrm{P}^0(\tau_k< k^2 \wedge \widehat{\tau}_0) \mathrm{P}^k(\widehat{\tau}_0>3t/2)}{\mathrm{P}^0(\widehat{\tau}_0>3t/2)} \\
& \ge & \frac{\mathrm{P}^0(\tau_k< k^2) \mathrm{P}^0( \tau_k < \widehat{\tau}_0) \mathrm{P}^k(\widehat{\tau}_0>3t/2)}{\mathrm{P}^0(\widehat{\tau}_0>3t/2)} \, .
\end{eqnarray*}
Now apply Lemma \ref{lem1}.\textbf{a}, \textbf{b} and \textbf{c} to the rightmost term to obtain a constant $h_2>0$ such that
for all $k$ and $t>2k^2$, the above term exceeds $h_2$. 
\hfill $\Box$

\bigskip

%%%%%%%%%%%%%%%%%%%%%%%%%%%%%%%%%%%%%%%%%%%%%%%%%%%%%%%%%%%%%%%%%%%%%%%%%%%%%%%%%
\section{Proof of Theorem \ref{thm1}} 
%%%%%%%%%%%%%%%%%%%%%%%%%%%%%%%%%%%%%%%%%%%%%%%%%%%%%%%%%%%%%%%%%%%%%%%%%%%%%%%%%
\label{sec:hybzones}
\setcounter{equation}{0}

Many arguments, notation and inequalities present in this section are borrowed from Cox and Durrett's paper \cite{CD}. It will be explicitly mentioned  ahead from which point our results take place. When a step is not explicitly justified the reader is invited to check \cite{CD}.

\medskip

Recall from Section \ref{sec:introd} the definitions of the random walk
$X^{x,t}$ that trace the origin of the value at $x$ at time $t$ and the dual process
$\left(\xi^{t}_{s}(A):\, 0\leq s \leq t \right)=\{X_{s}^{x,t}, x \in A:\, 0\leq s \leq t \}$, for $A\subset\ZZ$.

\medskip

We have that $r_t - l_t$ is bounded above by the number of ``inversions" at time $t$, i.e, the number of pairs $(x,y)$ so that $x<y$, $\eta_t(x)=0$ and $\eta_t(y)=1$. We denote this quantity by $B_t$. Then we should estimate $\mathrm{P}(B_t >M)$ for $M>0$. 

We start by observing that for an inversion at a pair of sites $x<y$ to happen at time $t$ the two dual coalescing random walks $X^{x,t}$ and $X^{y,t}$ cross each other without hitting and they end up with $X_{t}^{x,t}\geq 0>
X_{t}^{y,t}$. From here the clever idea in \cite{CD} was to realize that the event $\{B_t >M\}$, under the condition that no crossing of dual coalescing random walks happens after a dual time $K>0$ whose choice is independent of $t$, has probability of order $O(M^{-1})$. To state it precisely, let
\begin{eqnarray*}
  D_K(w,z,t) &=&  \left\{\xi^{t}_{K}(\ZZ) \cap [w,z]= \{w,z\}, X_{t-K}^{w,t-K} \geq 0 > X_{t-K}^{z,t-K}
  \right\},
\end{eqnarray*}\
\newline
for $t>K$ and $(w,z)\in\ZZ^2 \cap \{(u,v): u < v\}$ and put
\begin{eqnarray}
\label{eq:thm1-1}
  A_K(t) &=& \sum_{w<z} I_{D_K(w,z,t)}.
\end{eqnarray}
Then $A_K(t)$ counts the number of relevant changes of order that occur after time $K$. In \cite{CD} it was shown, under finite second moment assumption on the probability kernel $p(\cdot)$, that if $t>K$ then
$$
\mathrm{P}(B_t>M,A_K(t)=0) \le c \frac{K}{M}
$$
for a constant $c>0$ not depending on $t$, $K$ and $M$. Therefore
\begin{eqnarray*}
  \mathrm{P}(r_t-l_t>M) &\leq & \mathrm{P}(B_t>M)=\mathrm{P}(B_t>M ,A_K(t)>0)+ \mathrm{P}(B_t>M,A_K(t)=0) \\
    &\leq & \mathrm{P}(A_K(t)>0)+ \mathrm{P}(B_t>M,A_K(t)=0) \, .
\end{eqnarray*}
From this inequality we obtain Theorem \ref{thm1} as a consequence of the following result:

\bigskip

\begin{lemma} \label{lemm05}
For every $\epsilon >0 $, there exists $K>0$ so that $\mathrm{P}( A_K(t) >0)< \epsilon\,$  for all $t$
sufficiently large.
\end{lemma}

\bigskip

\begin{remark} Lemma \ref{lemm05} was proved in \cite{CD} under finite absolute third moment assumption on the probability kernel $p(\cdot)$. The proof we give only requires finite second moment.
\end{remark}

\bigskip

\no \textbf{Proof of Lemma \ref{lemm05}} 
As is usual, we fix $ \epsilon > 0$ but arbitrarily small.  The idea behind the proof is to consider the event $\{A_K(t)>0\}$ conditioned on $\xi^t_K(\ZZ)$ and them follow the dual coalescing random walks backward in time from dual time $K$. We then fix a positive integer $M$ that will be taken sufficiently large and we consider separatedly two groups of dual coalescing random walks with starting points taken from $\xi^t_K(\ZZ)$: those starting at $[-M\sqrt{t},M\sqrt{t}]$ and those starting outside this interval. The aim behind this partition of $\xi^t_K(\ZZ)$ is that for the second group the estimates can rely on the invariance principle while for the first one the relevant point is that the density of walks in $\xi^t_K(\ZZ) \cap [-M\sqrt{t},M\sqrt{t}]$ is small for large $K$.

For this discussion we introduce some notation: Let 
$$
D_1=\{ (x,y): x<y, \, \max(|x|,|y|)<M\sqrt{t} \}
$$ 
and
$$
D_2=\{ (x,y): x<y , \, \max(|x|,|y|) \ge M\sqrt{t} \} \, .
$$
By (\ref{eq:thm1-1}) we have that
\begin{equation}
\label{eq:thm1-2}
\mathrm{P}( A_K(t) >0) \le \mathrm{P} \left( \sum_{(w,z)\in D_1} I_{D_K(w,z,t)} > 0 \right) + \mathrm{P} \left( \sum_{(w,z)\in D_2} I_{D_K(w,z,t)} > 0 \right)
\end{equation}
We start with the first probability in the right hand side of (\ref{eq:thm1-2}). We are going to show that (given a fixed $M$) there exists a constant $c>0$ such that for $t > 2K$ sufficiently large 
\begin{equation}
\label{eq:thm1-21}
E \left[ \sum_{(w,z)\in D_1} I_{D_K(w,z,t)} \right] \le  \, \frac{\epsilon }{2} \, .
\end{equation}
By the Markov property
$$
I_1 := E\left( \sum_{(w,z)\in D_1} I_{D_K(w,z,t)} \, \Big| \, \xi^t_K(\ZZ) \right) = \sum_{(w,z)\in D_1} I_{\xi_K^{t}(\ZZ) \cap [w,z]=\{w,z\}} \mathrm{P}\left( X_{t-K}^{w,t-K}\geq 0>X_{t-K}^{z,t-K} \right).
$$
We will prove (\ref{eq:thm1-21}) by showing that with probability tending to one as $t$ tends to infinity, $I_1$ is small.

Recall from Section \ref{sec:introd} the definition of the family $(Y^x:x\in \ZZ)$ of independent random walks and the definition of $V^k(t)$. It follows directly from the definitions that $I_1$ is bounded above by
\begin{eqnarray*}
\lefteqn{ \!\!\!\!\!\!\!\!\!\!\!\!\!\!\!\!\!\!\!\!\!\!\!\!\!\!\!\!\!\!\!\!\!\!\!\!\!\!\!\!\!\!\!\!\!\!
\sum_{(w,z)\in D_1} I_{\xi_K^{t}(\ZZ) \cap [w,z]=\{w,z\}} \mathrm{P}\left( Y^w(t-K)>Y^z(t-K)\, , \, \tau_0^{Y^z-Y^w}>t-K \right) } \\
& & = \sum_{(w,z)\in D_1} I_{\xi_K^{t}(\ZZ) \cap [w,z]=\{w,z\}} \mathrm{P}\left( V^{z-w}(t-K) \right) \, .
\end{eqnarray*}
Denote $\mathrm{P} (V^k(t))$ by $f(k,t)$. By Lemma \ref{lem1}.\textbf{c}, there exists a constant $\gamma_3>0$ such that
\begin{equation}
\label{eq:thm1-3}
f(k,t) < \gamma_3 \, \frac{k}{\sqrt{t}} \, .
\end{equation}
Moreover, by Proposition \ref{lemm1} we can pick $k_0$ to satisfy 
\begin{equation}
\label{eq:thm1-4}
f(k,t) < \frac{\epsilon}{10M} \, \frac{k}{\sqrt{t}}
\end{equation}
for $t$ sufficiently large uniformly over $k \geq k_0$. 
%We choose $k_1>k_0$ so that with probability tending to one as $t$ tends to infinity,
%$$
%\sum_{(w,z)\in D_1} I_{\xi_K^{\ZZ,t}\cap [w,z]=\{w,z\}} I_{z-w>k_1}(z-w)/\sqrt{t-K} <\epsilon/H
%$$
%This can be done because $\xi_K^{\ZZ,t}$ is a stationary ergodic configuration not depending on $t$.
%Thus by Lemma \ref{lem1}
%\begin{eqnarray}
%\sum_{D_1\,:\, z-w>k_1} I_{\xi_K^{\ZZ,t}\cap [w,z]=\{w,z\}} f(z-w,t-K)&\leq&	 \sum_{D_1\,:\, k_1 \geq z-w>k_0} I_{\xi_K^{\ZZ,t}\cap [w,z]=\{w,z\}} c_3 \frac{z-w}{\sqrt{t-K}}\nonumber \\
%&\leq& \frac{\epsilon c_3}{H}. \label{A1}
%\end{eqnarray}
Now observe that for $t\geq 1$, deterministically,
$$
\sum_{(w,z)\in D_1} I_{\xi_K^{t}(\ZZ) \cap [w,z]=\{w,z\}}(z-w) \leq 2 M \sqrt{t},
$$ 
and therefore, for $t$ satisfying (\ref{eq:thm1-4}), we have that
\begin{eqnarray}
\label{eq:thm1-5}
\sum_{D_1\,: \, z-w>k_0} I_{\xi_K^{t}(\ZZ) \cap [w,z]=\{w,z\}} f(z-w,t-K)&\leq&	 \, \frac{\epsilon}{10M} \sum_{D_1\,: z-w>k_0} I_{\xi_K^{t}(\ZZ) \cap [w,z]=\{w,z\}} \, \frac{z-w}{\sqrt{t-K}}\nonumber\\
&\leq&  \frac{2 M \sqrt{t} \epsilon}{10M \sqrt{t-K}} < \epsilon / 4,
\end{eqnarray}
for $t$ large.
To deal with pairs $(w,z) \in D_1$ such that $z-w \leq k_0$ we apply (\ref{eq:thm1-3}) directly to obtain
\begin{equation}
\label{eq:thm1-6}
\sum_{D_1\,:\, z-w \leq k_0} I_{\xi_K^{t}(\ZZ) \cap [w,z]=\{w,z\}} f(z-w,t-K) \leq \frac{\gamma_3k_0}{\sqrt{t-K}} |\xi^t_K \cap (-M \sqrt t , M \sqrt t )| \, .
\end{equation}

Now the random variable $|\xi^t_K \cap (-M \sqrt t , M \sqrt t )|$ satisfies 
$\lim_{t \rightarrow \infty} \frac{|\xi^t_K \cap (-M \sqrt t , M \sqrt t )|}{2M \sqrt t}  = \delta(K)$ in
probability, where $\delta (K)$ is the density of our system of coalescing random walks at time $K$. We simply now use the well known fact that $\delta(K)$ goes to zero as $K$ becomes large to conclude that with probability tending to one as $t$ tends to infinity
$$
\sum_{D_1\,:\, z-w \leq k_0} I_{\xi_K^{t}(\ZZ) \cap [w,z]=\{w,z\}} f(z-w,t-K) \leq 2\delta(K)2M \sqrt t \frac{k_0\gamma_3}{\sqrt {t-K}} <
\epsilon/2,
$$
provided $K$ is fixed large. This implies that with this choice of $K$, $\mathrm{E} \left( \sum_{(w,z)\in D_1} I_{D_K(w,z,t)} > 0 \right) < \epsilon /2.$
%To bound the expectation of right hand side summation in the previous inequality we use the fact that the %density of $\xi_K^{t}(\ZZ)$ tends to zero as $K$ tends to infinity. We have that
%\begin{eqnarray}
%\label{eq:thm1-7}
%\lefteqn{
%E\left[\sum_{(w,z)\in D_1} I_{\xi_K^{t}(\ZZ) \cap [w,z]=\{w,z\}}(z-w) I_{z-w \leq k_0} \right] } \nn \\
%&\leq & k_0 E\left[\sum_{(w,z)\in D_1} I_{\xi_K^{t}(\ZZ) \cap [w,z]=\{w,z\}} I_{z-w \leq k_0} \right] \nn \\
%& \le & k_0 \frac{2M\sqrt{t}}{k_0} \, \mathrm{P} \left( \{0,1\} \subset \xi_K^{t}(\ZZ) \right) \nn \\
%&\le & 2M\sqrt{t} \, \mathrm{P} \left( 0 \in \xi_K^{t}(\ZZ) \right)^2 \le c \, \frac{M\sqrt{t}}{K}
%\end{eqnarray}
%for some constant $c>0$. In the previous expression the third inequality comes from a more general FKG type property obtained by Arratia, see Lemma 1 in \cite{A2} and also Lemma 2.0.8 is \cite{Sun} for a discussion on Arratia's result as well as a generalization for discrete time, and for the fourth inequality see Lemma 2.0.7 in \cite{Sun}.
%Therefore, putting together (\ref{eq:thm1-5}), (\ref{eq:thm1-6}) and (\ref{eq:thm1-7}) we have shown that
%$$
%E[I_1] \leq (2+c\gamma_3)\frac{M \sqrt{t}}{K \sqrt{t-K}} \le \sqrt{2}(2+c\gamma_3) \, \frac{M}{K}.
%$$
%for every $t>2K$ sufficiently large, which gives (\ref{eq:thm1-21}). 

\medskip
Now we deal with the second probability in the right hand side of (\ref{eq:thm1-2}). We are going to show that there exists a constant $c>0$ such that 
\begin{equation}
\label{eq:thm1-22}
\mathrm{P} \left( \sum_{(w,z)\in D_2} I_{D_K(w,z,t)} > 0 \right) \le c \, E\left( \left(\frac{Y^0(2t)}{\sqrt{t}} -(M+1)\right)_+ \right)
\end{equation}
for every $K$ and $t$ sufficiently large. 

We fix integer $m$ such that $m \geq M$ or $m < -M-1$. Then write $m \sqrt{t}\leq x_0 < x_1...< x_R \leq (m+1)\sqrt{t}$ for
$\xi_K^{t}(\ZZ) \cap [m\sqrt{t}, (m+1)\sqrt{t}]$. Let $\Gamma(m,t)$ be the event that there exists consecutive coalescing random walks at $x_i, x_{i+1}$ at dual time $K$ so that
\begin{enumerate}
    \item [(i)] The random walks have not coalesced by dual time $t$
    \item [(ii)] $x_i, x_{i+1}\in [m \sqrt{t}, (m+1)\sqrt{t}]$
    \item [(iii)] If $m>0$ then at least one of the random walks is in $(-\infty, 0)$ at dual time $t$ and if $m<0$ then at least one of the random walks is in $(0,+\infty)$ at dual time $t$.
\end{enumerate}
Put for $m > 0$
$$
\tilde{\Gamma}(m,t) = \Gamma(m,t) \cup \{ X^{x_0, t-K}_{t-K} \textrm{ or } X^{x_R, t-K}_{t-K} \textrm{ is in } (-\infty, 0) \textrm{ at time t}\},
$$
and for $m < 0$
$$
\tilde{\Gamma}(m,t) = \Gamma(m,t) \cup \{ X^{x_0, t-K}_{t-K} \textrm{ or } X^{x_R, t-K}_{t-K} \textrm{ is in } ( 0, \infty) \textrm{ at time t}\},
$$
and observe that
$$
\left\{ \sum_{(w,z)\in D_2} I_{D_K(w,z,t)} > 0 \right\} \subset \bigcup _{m \in \ZZ-(-M,M] } \tilde{\Gamma}(m,t)
$$
Then by symmetry we have that the probability in (\ref{eq:thm1-22}) is bounded above by
$$
2 \mathrm{P}\left( \bigcup _{m = M}^{+\oo} \tilde{\Gamma}(m,t) \right) \le 2 \sum_{m=M}^{+\oo} \mathrm{P}\left( \tilde{\Gamma}(m,t) \right).
$$

By the reflection principle the probability of the second event in the definition of $\tilde{\Gamma}(m,t)$ is bounded above by 
\begin{equation}
\label{eq:thm1-8}
2 \mathrm{P}\left( \sup_{0 \leq s \leq t} Y^0(s) \geq m \sqrt{t}\right)
\le 4 \mathrm{P}\left( |Y^0(2t)| \geq m \sqrt{t}\right),
\end{equation}
We will see that a similar bound applies to the probability of $\Gamma (m,t)$. Conditioning on the event that $\xi_K^{t}(\ZZ) \cap [m \sqrt{t}, (m+1)\sqrt{t}]=\{x_0, x_1,...,x_R\}$, where $R$ is fixed integer, we have by Lemma \ref{lem1}.\textbf{c} that the probability that the random walks $(X^{x_i,t-K}(s)\,:\, 0\leq s \leq t-K)$ and
$(X^{x_{i+1},t-K}(s)\,:\, 0\leq s \leq t-K)$ have not coalesced by dual time $t-K$ is bounded by $
\gamma_3 (x_{i+1}-x_i)/\sqrt{t}$. Moreover, by Proposition \ref{lemm2}, given that these two random walks do not coalesce, the probability that either random walk ends to the left of 0 is less than 
$$
c_0 P\left(  |Y^0(2t)| \geq  m \sqrt{t}\right) \, .
$$ 
Hence the conditional probability of $\Gamma (m,t)$ given $\xi_K^{t}(\ZZ) \cap [m \sqrt{t}, (m+1)\sqrt{t}]=\{x_0, x_1,...,x_R\}$ is bounded above by
\begin{equation}
\label{eq:thm1-9}
c_0 \gamma_3 \left( \sum_{i=1}^{R} \frac{x_i-x_{i-1}}{\sqrt{t}} \right) \mathrm{P}\left( |Y^0(2t)| \geq  m \sqrt{t}\right) \le c_0 \gamma_3  \mathrm{P}\left( |Y^0(2t)| \geq  m \sqrt{t}\right).
\end{equation}
Thus we apply (\ref{eq:thm1-8}) and (\ref{eq:thm1-9}) to get
\begin{eqnarray}
\sum_{m=M}^{+\oo} \mathrm{P}\left( \tilde{\Gamma}(m,t) \right) & \leq &  
(c_0 \gamma_3+4) \sum_{m=M}^{+\oo} \mathrm{P}\left( |Y^0(2t)| \geq  m \sqrt{t}\right) \nn \\
& \le & (c_0 \gamma_3+4) E\left(
\left(\frac{Y^0(2t)}{ \sqrt{t}} -(M+1)\right)_+ \right), \nn
\end{eqnarray}
and we have obtained (\ref{eq:thm1-22}) with $c=2 (c_0 \gamma_3+4)$.

\medskip
Since $\left({Y^0(2s)}/{\sqrt{s}}\right)_{s\geq 2K}$ is a uniformly integrable family of random variables (the $L_2$ norms are uniformly bounded), by the invariance principle we can chose $M$
sufficiently large such that
$$
E\left( \left(\frac{Y^0(2t)}{\sqrt{t}} -(M+1)\right)_+ \right) < \frac{\epsilon}{2 c},
$$ 
for all large $t$.

Thus with this choice of $M$ we have (for any choice of $K$, that for $t$ large
$$
\sum_{D_2\,:\, z-w \leq k_0} I_{\xi_K^{t}(\ZZ) \cap [w,z]=\{w,z\}} f(z-w,t-K) \ < 
\epsilon/2,
$$
We now fix $K$ sufficiently large to ensure that we have a similar inequality (for this value of $M$ for the summation over $D_1$.  Thus 
$$
\mathrm{P}( A_K(t) >0) \le  \epsilon \, .
$$
\hfill $\square$

\bigskip

%%%%%%%%%%%%%%%%%%%%%%%%%%%%%%%%%%%%%%%%%%%%%%%%%%%%%%%%%%%%%%%%%%%%%%%%%%%%%%%%%%%%
\section{Proof of Theorem \ref{thm2}} 
%%%%%%%%%%%%%%%%%%%%%%%%%%%%%%%%%%%%%%%%%%%%%%%%%%%%%%%%%%%%%%%%%%%%%%%%%%%%%%%%%%%%
\label{sec:theorem2}
\setcounter{equation}{0}

To prove that $(r_t - l_t : t \ge 0)$ is not tight we have to show that there exists a constant $\theta > 0$ such that for a strictly increasing sequence of positive integers $(M_k)_{k \ge 1}$ and a sequence $(t_k)_{k \ge 1}$ of positive real numbers
\begin{equation}
\label{eq:notight1}
\limsup_k \mathrm{P} ( r_{t_k} - l_{t_k} \ge M_k ) \ge \theta \, .
\end{equation}
In this section we will perform the claculations under the assumption that the kernal $p( \cdot)$ is symmetric.  This is to simplify the calculations and is not necessary, we leave to the reader the messy but straightforward extra steps
to arrive at the general case.

Now recall from Section \ref{sec:introd} the definition of the families of random walks $(X^{x,t}: x\in \ZZ, \, t\ge 0)$ and $(Y^x:x\in \ZZ)$. By duality, we have that 
\begin{equation}
\label{eq:notight2}
\mathrm{P} ( r_{t_k} - l_{t_k} \ge M_k ) = \mathrm{P}( \textrm{There exists } x < x + M_k \le y \textrm{ such that } X^{x,t_k}_{t_k} > 0 \ge X^{y,t_k}_{t_k}), 
\end{equation}
and indeed we obtain (\ref{eq:notight1}) by showing that for appropriate sequences $(M_k)_{k \ge 1}$ and $(t_k)_{k \ge 1}$, we have
$$
\limsup_k \mathrm{P}( X^{M_k,t_k}_{t_k} > 0 \ge X^{3M_k,t_k}_{t_k} ) \ge \theta ,
$$
which is stronger than (\ref{eq:notight1}) since the last probability is clearly bounded below by the right hand side in (\ref{eq:notight2}). Now
\begin{equation}
\label{eq:notight3}
\mathrm{P}( X^{M_k,t_k}_{t_k} > 0 \ge X^{3M_k,t_k}_{t_k} ) = \mathrm{P}( Y^{M_k}_{t_k} > 0 \ge Y^{3M_k}_{t_k}, \, \tau_0^{Y^{3M_k}-Y^{M_k}} > t_k ).
\end{equation}
Write
$$
Y^{M_k} = Z_{k,1}^{\prime} + Z_{k,1}^{\prime \prime} + M_k
\qquad \textrm{and} \qquad
Y^{3M_k} = Z_{k,2}^{\prime} + Z_{k,2}^{\prime \prime} + 3 M_k
$$
where, $Z_{k,i}^{\prime}$, and $Z_{k,i}^{\prime \prime}$, $i=1,2$, are independent symmetric random walks starting at $0$ with transition kernels given respectively by
$$
p^{\prime} (x) = p(x) \textbf{1}_{|x| \le 4 M_k} \qquad \textrm{and} \qquad 
p^{\prime \prime} (x) = p(x) \textbf{1}_{|x| > 4 M_k}.
$$
Let $F_{k,1}$ be the event that during the time interval $[0,t_k]$ the random walk $Z_{k,1}^{\prime \prime}$ makes a unique jump of size greater or equal to $4M_k$
and that this jump occurs to the right and $F_{k,2}$ the event that during the time interval $[0,t_k]$ the random walk $Z_{k,2}^{\prime \prime}$ makes a unique jump of size greater or equal to $4M_k$ and that this jump occurs to the left. The right hand side of (\ref{eq:notight3}) 
is bounded below by
\begin{eqnarray}
\label{eq:notight4}
\lefteqn{ \!\!\!\!\!\!\!\!\!\!\!\!\!\!\!\!\!\!\!\!\!\!\!\!
\mathrm{P}\left( F_{k,1} \cap F_{k,2}, \, \sup_{ 0 \leq s \leq t_k} |Z_{k,1}^{\prime}(s)|< M_k,\, \sup_{ 0 \leq s \leq t_k} |Z_{k,2}^{\prime}(s)|< M_k \right) } \nn \\
& & =  \mathrm{P}(F_{k,1})^2 \, \mathrm{P} \left( \sup_{ 0 \leq s \leq t_k} |Z_{k,1}^{\prime}(s)|< M_k \right)^2 \, .
\end{eqnarray}
where the equality comes from independence and symmetry. From this point we fix the sequences $(M_k)_{k \ge 1}$ and $(t_k)_{k \ge 1}$. We put for every $k \ge 1$
$$
M_k = 2^k
\qquad \textrm{and} \qquad
t_k = C \left( \sum_{x \ge 2^{k+2}} p(x) \right)^{-1} \, ,
$$
where $C$ is a constant that are going to be fixed later. With these definitions it is simple, by calculating the jump rates, to obtain that for every $k \ge 1$
$$
\mathrm{P} (F_{k,1}) = C e^{-2C} \, .
$$
Therefore to finish the Proof of Theorem \ref{thm2}, by reversing the inequality in the expression inside the last probability in (\ref{eq:notight4}), we only have to show that there exists a constant $0 < \vartheta < 1$ such that
\begin{equation}
\label{eq:notight5}
\limsup _k \mathrm{P} \left( \sup_{ 0 \leq s \leq t_k} |Z_{k,1}^{\prime}(s)| \ge 2^{k} \right) < \vartheta \, ,
\end{equation}
for an appropriate choice of $C$.

Using the symmetry and the Strong Markov property, we have
\begin{equation}
\label{eq:notight6}
\mathrm{P} \left( \sup_{ 0 \leq s \leq t_k} |Z_{k,1}^{\prime}(s)| \ge 2^{k} \right) \le 4 P(Z_{k,1}^{\prime}(t_k)> 2^{k})
\le 2^{2(1-k)}  E[Z_{k,1}^{\prime}(t_k)^2] \, .
\end{equation}
Let $R$ be an integer valued random variable with distribution given by $p(\cdot)$. Since 
$$
Z_{k,1}^{\prime}(t)^2 - t \sum_{x=-2^{k+2}}^{2^{k+2}} x^2 p(x)
$$
is a martingale, we have that 
\begin{eqnarray}
E[Z_{k,1}^{\prime}(t_k)^2] &=& t_k \sum_{x=-2^{k+2}}^{2^{k+2}} x^2 p(x) \nn \\
    & \leq& 2t_k \sum_{i=1}^{k+2} 2^{2i} \, \mathrm{P}(|R| \in [2^{i-1}, \, 2^i]) \nn \\
    &=& 2 C \frac{\sum_{i=1}^{k+2} 2^{2i} \, \mathrm{P}(|R| \in [2^{i-1}, \, 2^i])}{\sum_{|x|>2^{k+2}}p(x)}. \nn
\end{eqnarray}
We claim that there exists a constant $C^\prime>0$ such that 
$$
\sup_k \frac{\sum_{i=1}^{k} 2^{2i} \, \mathrm{P}(|R| \in [2^{i-1}, \, 2^i])}{2^{2k} \sum_{|x|>2^k}p(x)} \le C^\prime.
$$
Now chose $C>0$ such that $2^6 C C^{\prime} < 1$ and take $\vartheta = 2^6 C C^{\prime}$ in (\ref{eq:notight5}) to conclude from (\ref{eq:notight6}).

\medskip

It remains to prove the claim. Here we recall our basic hypotheses that there exists $\epsilon \in (0,1)$ such that $E[|R|^{2-\epsilon}] = \sum |x|^{2-\epsilon} p(x) = +\oo$. Clearly this is equivalent to 
$$
I_{\epsilon}=\sum_{k>0} 2^{(2-\epsilon)k} \mathrm{P}\left(|R| \in [2^{k-1},
2^{k}]\right) = +\oo \, .
$$
So reducing $\epsilon$ slightly if necessary, we have $\limsup_{k\rightarrow \infty} 2^{(2-\epsilon)k} \mathrm{P}(|R| \in [2^{k-1}, \, 2^k])=\infty$. Let as denote by $V_k=2^{2k} \mathrm{P}(|R| \in [2^{k-1}, \, 2^{k}])$, then, there exists a subsequence $k_i \rightarrow \infty$ for which we have $ V_{k_{i}} 2^{-k_i \epsilon } = \sup_{j \leq k_i} V_j
2^{-j \epsilon }.$ Thus
\begin{eqnarray*}
  \frac{{\sum_{j=1}^{k_i-1} 2^{2j} \, \mathrm{P}(|R| \in [2^{j-1}, \, 2^j])}}{{2^{2k_{i}} \mathrm{P}(|R| \in [2^{k_i-1}, \, 2^{k_{i}}])}} &=&   \frac{{\sum_{j=1}^{k_i-1} V_j}}{ 2^{2k_i} \mathrm{P}(|R| \in [2^{k_i-1}, \, 2^{k_{i}}])} \\
    &\leq& \frac{2^{-k_{i} \epsilon} V_{k_{i}}{\sum_{j=1}^{k_i-1}  2^{j \epsilon}}}{ {2^{2k_i} \mathrm{P}(|R| \in [2^{k_i-1}, \, 2^{k_{i}}])} }  \\
    &\leq & \sum_{j=1}^{k_i-1}  2^{-(k_{i}-j) \epsilon} < C^{\prime} \, ,
\end{eqnarray*}
for some $C^{\prime} = C^{\prime}(\epsilon)$. The proof now is complete. $\Box$

\appendix
%%%%%%%%%%%%%%%%%%%%%%%%%%%%%%%%%%%%%%%%%%%%%%%%%%%%%%%%%%%%%%%%%%%%%%%%%%%%%%%%%%%%%%
\section{Basic Results on Random Walks}
%%%%%%%%%%%%%%%%%%%%%%%%%%%%%%%%%%%%%%%%%%%%%%%%%%%%%%%%%%%%%%%%%%%%%%%%%%%%%%%%%%%%%%
\label{sec:basres}
\setcounter{equation}{0}

In this section we will assemble some well known properties of random walks, most of which come from the local central limit theorem for sums of i.i.d. random variables with second moment.  The standard reference is \cite{Spi}. 

Analogously to the notation introduced in Section \ref{sec:introd} we define $\tau_A^S=\inf \{n \ge 0\, :\; S_n \in A\}$ and $\widehat{\tau}_A^S=\inf \{n > 0\, :\; S_n \in A\}$ for a discrete time stochastic process $(S_n)_{n\geq 0}$ and $A \subset \ZZ$.
%show some very interesting properties which relate to the stopping time,
%the comparison between two stopping times and its asymptotic behavior.\\

\begin{lemma} 
\label{lem1}
Let $(X_n)_{n\geq 0}$ be a one-dimensional symmetric homogeneous translation invariant discrete time random walk whose transition function $p(x)$ has finite second moment, then
\begin{enumerate}
\item[\textbf{(a)}] For every $x \neq 0$ we have that $\mathrm{P}^0(\tau_x < \widehat{\tau}_0)\approx 1/|x|.$
\item[\textbf{(b)}] $\mathrm{P}^0(\widehat{\tau}_0 \geq n^2 ) \approx 1/n$.
\item[\textbf{(c)}] Uniformly in $x\neq 0$, $\mathrm{P}^x(\tau_0 \geq n )\approx (|x| /\sqrt{n})\wedge1$.
\item[\textbf{(d)}] There exists $\gamma_5>0$ so that for all $0 \leq x \leq k$ we have $\mathrm{P}^x(\tau_k<\tau_0)\leq \gamma_5 x/k.$
\end{enumerate}
In the case of continuous time processes, the results remain the same.
\end{lemma}

\no \textbf{Proof:} 

\smallskip

\no \textbf{(a)} Let $g_0$ be the green function for the process killed at 0, defined as follows:
$$ g_0(x,y) =\left\{%
\begin{array}{ll}
    \sum_{n=0}^\infty \mathrm{P}^x(X_n=y , \tau_0>n)  & \hbox{ if $x$ and $y \in \ZZ-\{0\}$} \\
    0 & \hbox{otherwise.} \\
\end{array}%
\right. 
$$
Let, for $m \geq 1$, $ \tau_x^m=\inf \{r>\tau_x^{m-1}\, \,\mbox{so that} \,\, X_r=x \}$, with
$\tau_x^0=0$ and
$ \tau_x^1 = \tau_x$.\\
Then for $x\neq 0$:
\begin{eqnarray*}
  g_0(x,x) &=& \sum_{n=0}^\infty \mathrm{P}^x(X_n=x , \tau_0>n) = 1+\sum_{n=1}^\infty \mathrm{P}^x(X_n=x , \tau_0>n) \\
   &=& 1+\sum_{n=1}^{\infty} {\sum_{m=1}^n \mathrm{P}^x(\tau_x^m=n ,\tau_0>n)} = 1+\sum_{m=1}^{\infty}  \mathrm{P}^x(\tau_x^m <\tau_0)  \\
    &=&1+\sum_{m=1}^{\infty}  \left(\mathrm{P}^x(\widehat{\tau}_x <\tau_0)\right)^m = \frac{1}{\mathrm{P}^x(\tau_0<\widehat{\tau}_x)}. \\
\end{eqnarray*} \
\newline
We now observe that $g_0(x,x)=2a(x)$, (see P29.4 in \cite{Spi}), where 
$$
a(x)= \sum_{i=0}^\infty \left(\mathrm{P}^0(X_i=0)-\mathrm{P}^x(X_i=0)\right),
$$ 
is the potential kernel function. Moreover, by  P28.4 in \cite{Spi}, $a(x)/x \rightarrow \sum_x x^2 p(x)$ as $x \rightarrow \infty$, i.e., $a(x)\approx x$. Therefore
$$
\mathrm{P}^x(\tau_0<\widehat{\tau}_x)=\frac{1}{2a(x)} \approx \frac{1}{|x|} \, .
$$

\medskip

\no \textbf{(b)} The result follows from P.32.3 in \cite{Spi} which state that
$$
\lim_{n\longrightarrow \infty} \sqrt{n}\mathrm{P}^0(\tau_0 >n)= \sqrt{\frac{2}{\pi}}
\left( \sum_x x^2 p(x) \right)^{\frac{1}{2}} .
$$
\no \textbf{(c)} We deduce from \textbf{b} that there exists constants $0< \gamma_2< \gamma_1$ , such that 
$$
\frac{\gamma_2}{n} \leq \mathrm{P}^0(\tau_0 \geq n^2 )\leq \frac{\gamma_1}{n},
$$ 
for each positive integer $n$. We also have that
$$  
\mathrm{P}^0(\widehat{\tau}_0>n^2) \geq \mathrm{P}^0(\tau_x <\widehat{\tau}_0)\mathrm{P}^x(\tau_0>n^2)
$$
and therefore we obtain the upper bound
$$
\mathrm{P}^x(\tau_0>n^2) \le \frac{\mathrm{P}^0(\widehat{\tau}_0> n^2)}{\mathrm{P}^0(\tau_x<\widehat{\tau}_0)}
\le \frac{c_1}{n \mathrm{P}^0(\tau_x<\widehat{\tau}_0)} \le \gamma_3 \frac{|x|}{n} \, .
$$
So $\mathrm{P}^x(\tau_0 > n^2) \leq \gamma_3 ( \frac{|x|}{n}\wedge 1) $ for a constant $\gamma_3>0$ not depending on $x$ or $n$.

For the lower bound, we use a last passage time at $0$ decomposition. First we write
\begin{eqnarray*}
 \mathrm{P}^x(\tau_0 \leq n^2) &=& \sum_{k=0}^{n^2} \mathrm{P}_k(x,0) \mathrm{P}^0(\widehat{\tau}_0 >n^2-k) \;\;\ \mbox{ for } x\neq 0\\
 1 &=&  \sum_{k=0}^{n^2} \mathrm{P}_k(0,0) \mathrm{P}^0(\widehat{\tau}_0 >n^2-k).
\end{eqnarray*}
Thus
\begin{eqnarray*}
   \mathrm{P}^x(\tau_0 > n^2) &=& \sum_{k=0}^{n^2} \left(\mathrm{P}_k(0,0)-\mathrm{P}_k(x,0)\right) \mathrm{P}^0(\widehat{\tau}_0 >n^2-k) \\
   &\geq& \mathrm{P}^0(\widehat{\tau}_0>n^2) \sum_{k=0}^{n^2} \left(\mathrm{P}_k(0,0)-\mathrm{P}_k(x,0)\right) \\
   &\geq& \gamma_4 \frac{|x|}{n}\wedge 1.
\end{eqnarray*}

\medskip

\no \textbf{(d)}
By the Markov Property
$$
 \mathrm{P}^0(\tau_k<\widehat{\tau}_0) \geq  \mathrm{P}^0(\tau_x<\widehat{\tau}_0) \mathrm{P}^x(\tau_k<\tau_0).
$$
Then $\mathrm{P}^x(\tau_k<\tau_0)\leq \mathrm{P}^0(\tau_k<\widehat{\tau}_0)/\mathrm{P}^0(\tau_x<\widehat{\tau}_0)$, and, by applying \textbf{a}, we find the required result.

\medskip

The extension of these results continuous processes is either
automatic, as part in \textbf{a} and \textbf{d}, or trivial.   $\square$

\end{document}